\documentclass[12pt]{article}

\usepackage{fullpage}
\usepackage{amssymb}
\usepackage{amsmath}
\usepackage[usenames]{color}
\usepackage{graphicx}

\usepackage{wrapfig}

\usepackage{multirow}
\usepackage{booktabs}

\usepackage[colorlinks=true,
linkcolor=webgreen,
filecolor=webbrown,
citecolor=webgreen]{hyperref}

\definecolor{webgreen}{rgb}{0,.5,0}
\definecolor{webbrown}{rgb}{.6,0,0}
\definecolor{red}{rgb}{1,0,0}

\usepackage{color}

\begin{document}
	
\title{Motzkin Numbers: an Operational Point of View}

\author{M. Artioli\\ ENEA - Bologna Research Center\\
	 Via Via Martiri di Monte Sole,  4, 40129, Bologna, Italy\\
	 marcello.artioli@enea.it 
	\and  G. Dattoli \\
	 ENEA - Frascati Research Center,\\
	  Via Enrico Fermi 45, 00044, Frascati, Rome, Italy\\
	 giuseppe.dattoli@enea.it 
	\and  S. Licciardi\\
	Dep. of Mathematics and Computer Science,\\
	 University of Catania, Viale A. Doria 6, 95125, Catania, Italy \\
	 ENEA - Frascati Research Center,\\
	  Via Enrico Fermi 45, 00044, Frascati, Rome, Italy \\
	  silvia.licciardi@dmi.unict.it 
	\and S. Pagnutti\\
	ENEA - Bologna Research Center,\\
	 Via Via Martiri di Monte Sole,  4, 40129, Bologna, Italy \\
	simonetta.pagnutti@enea.it}		
		
\maketitle		
\begin{abstract}
The Motzkin numbers can be derived as coefficients of hybrid polynomials. Such an identification allows the derivation of new identities for this family of numbers and offers a tool to investigate previously unnoticed links with the theory of special functions and with the relevant treatment in terms of operational means. The use of umbral methods opens new directions for further developments and generalizations, which leads, e.g., to the identification of new Motzkin associated forms.
\end{abstract}
		
%

\section{Introduction}
The telephone numbers $(T_{n})$, also called convolution numbers, provide a very well known example of link between special numbers and special polynomials. The $(T_{n})$ can be expressed in terms of Hermite polynomials coefficients $ (h_{s}) $  \cite{Riordan}. Two of the present authors (M.A and G.D.) have recently pointed out in ref. \cite{Artioli} that the Padovan and Perrin numbers \cite{Perrin,Padovan} can be recognized to be associated with particular values of two variable Legendre polynomials \cite{Dattoli}.\newline

Weinstein has discussed in \cite{Weisstein} the connection between Motzkin numbers and a family of hybrid polynomials, and Blasiak et al. and Dattoli et al. have studied, in \cite{Blasiak,Lorenzutta}, the relevant properties of Motzkin numbers.\newline

\noindent The hybrid polynomials are indeed defined as \cite{Lorenzutta}

\begin{equation}\label{Kn}
P_{n}^{(q)}(x,y)=n!\sum_{r=0}^{\lfloor\frac{n}{2}\rfloor }\dfrac{x^{n-2r}y^{r}}{(n-2r)!r!(r+q)!},
\end{equation}
and the relevant generating function reads

\begin{equation}\label{genKn}
\sum_{n=0}^{\infty}\dfrac{t^{n}}{n!}P_{n}^{(q)}(x,y)=\dfrac{I_{q}(2 \sqrt{y}\; t)}{(\sqrt{y}\;t)^{q}}e^{xt},
\end{equation}
where $I_{q}(x)$  is the modified Bessel function of the first kind of order $ q $.\newline

\noindent Within the present framework, the Motzkin numbers sequence can be specified as \cite{Blasiak}              

\begin{equation}\begin{split}\label{key}
& m_{n} = P_{n}^{(1)}(1,1)=\sum_{s=0}^{n}m_{n,s},\\
& m_{n,s}=\binom{n}{s}\;f_{s},\\
& f_{s} = \dfrac{s!}{\Gamma\left( \dfrac{s}{2}+2\right) \Gamma\left( \dfrac{s}{2}+1\right) }\left| \cos\left( s\dfrac{\pi}{2}\right) \right|, 
\end{split}\end{equation}
 where the coefficients $m_{n,s}$ can be represented as the triangle reported in the following table, in which $m_{n,2}$ corresponds, in OEIS, to the sequence $A000217$, $m_{n,4}$ to $A034827$, $m_{n,6}$ to $A000910$ and so on.\\

\begin{table}[h]
	\centering
\begin{tabular}{||c | c|c|c|c|c|c|c|c|c||c||} \hline
	\multicolumn{10}{||c||}{\centering $m_{n,s}$ \textbf{coefficients}}& $ m_{n}\; \textbf{Motzkin} $ \\ \hline
	\toprule
	 \multirow{2}{*}{\textbf{Parameter}}  & \multicolumn{9}{c||}{\textbf{s}}&\multirow{2}{*}{$\sum_{s=0}^{n}m_{n,s}$}\\ 
	 \cmidrule(lr){3-10}
	&  & \textbf{0} & \textbf{1} & \textbf{2 } & \textbf{3} & \textbf{4} & \textbf{5} & \textbf{6} & \textbf{7}&\\  
	\midrule	
	 \multirow{9}{*}{\textbf{n}}  & \textbf{0} & 1 &  &  &  & & & & & \textbf{1} \\ 
	& \textbf{1}  & 1 & 0 &  &  & & &  & & \textbf{1} \\ 
	& \textbf{2} & 1 & 0 & 1 & & & &  & & \textbf{2} \\ 
	 & \textbf{3} & 1 & 0 & 3 & 0 &  &  & & & \textbf{4}\\ 
	 & \textbf{4} & 1 & 0 & 6 & 0 & 2 &  & & & \textbf{9}\\ 
	& \textbf{5} & 1 & 0 & 10 & 0 & 10 & 0  &  &  & \textbf{21}\\ 
	& \textbf{6} & 1 & 0 & 15 & 0 & 30 & 0 & 5 &  & \textbf{51}\\ 
	& \textbf{7} & 1 & 0 & 21 & 0 & 70 & 0 & 35 & 0  & \textbf{127}\\ 
	& \dots &  \dots & \dots & \dots    &\dots  &\dots  &\dots    &\dots  &\dots  & \dots\\ \hline
	\bottomrule
\end{tabular}
\caption{Motzkin Numbers and their  Coefficients.}
\label{table1}
\end{table}

According to eq. \eqref{genKn}, the Motzkin numbers can also be defined as the coefficients of the following series expansion

\begin{equation}\label{genMn}
\sum_{n=0}^{\infty}\dfrac{t^{n}}{n!}m_{n}=\dfrac{I_{1}(2  t)}{t}e^{t}.
\end{equation}

In the following we will show how some progresses in the study of the relevant properties can be done by the use of a formalism of umbral nature.\newline

\section{Motzkin Numbers and Umbral Calculus}
 
 In order to simplify most of the algebra associated with the study of the properties of the Motzkin numbers and to get new relevant identities, we introduce a formalism successfully exploited elsewhere \cite{Borel} based on methods of umbral nature \cite{Roman}. \newline
 To this aim we note that the function

\begin{equation}\label{Tric}
C_{q}(x)=\dfrac{I_{q}(2 \sqrt{x})}{(\sqrt{x})^{q}}=\sum_{r=0}^{	\infty}\dfrac{x^{r}}{r!(q+r)!}
\end{equation}
can be cast in the form

\begin{equation}\label{TricOp}
C_{q}(x)=\hat{c}^{q} \circ e^{\hat{c}x},
\end{equation}
where $\hat{c}$  is an umbral operator defined according to

\begin{equation}\label{actOp}
\hat{c}^{\mu}=\dfrac{1}{\Gamma(\mu+1)},
\end{equation}
with $\mu$ not necessarily integer and real.\newline
We define the following composition rule 

\begin{equation}\label{opop}
\hat{c}^{\mu} \circ \hat{c}^{\nu}=\hat{c}^{\mu+\nu}
\end{equation}
and we let $\hat{C}=\{\hat{c}^{\alpha},\alpha\in \mathbb{C}\}$ denote the set of $\hat{c}$-operators. Then,  the pair $(\hat{C},\circ)$ satisfying the Abelian-group property. The mathematical foundations of the theory of $\hat{c}$-operators can be traced back to those underlying the Borel transform and have been carefully discussed in ref. \cite{Roman}.\newline

The use of this formalism allows to restyle the hybrid polynomials in the form

\begin{equation}\label{KnOp}
P_{n}^{(q)}(x,y)=\hat{c}^{q} \circ H_{n}(x,\hat{c}\;y),
\end{equation}
where

\begin{equation}\label{Herm}
H_{n}(x,y)=n!\sum_{r=0}^{\lfloor\frac{n}{2}\rfloor }\dfrac{x^{n-2r}y^{r}}{(n-2r)!r!}
\end{equation}
are the two variable Hermite-Kamp\'e de F\'eri\'et  polynomials of order $ 2 $.\newline
We can accordingly use the wealth of properties of this family of polynomials to derive further and new relations regarding those of the Motzkin numbers family.\newline

\noindent By recalling indeed the generating function \cite{Lorenzutta}

\begin{equation}\label{genHermnl}
\sum_{n=0}^{\infty}\dfrac{t^{n}}{n!}H_{n+l}(x,y)=H_{l}(x+2yt,y)e^{xt+yt^{2}},
\end{equation}
we find

\begin{equation}\label{genmnl}
\sum_{n=0}^{\infty}\dfrac{t^{n}}{n!}m_{n+l}=\hat{c} \circ H_{l}(1+2\hat{c}t,\hat{c})e^{t+\hat{c}t^{2}},
\end{equation}
which, after using eqs. \eqref{opop}, \eqref{Herm}, \eqref{TricOp}, finally yields

\begin{equation}\begin{split}\label{mu}
& \sum_{n=0}^{\infty}\dfrac{t^{n}}{n!}m_{n+l}=\mu_{l}(t)\;e^{t},\\
& \mu_{l}(t)=l!\sum_{r=0}^{\lfloor\frac{l}{2}\rfloor }\dfrac{1}{r!}\sum_{s=0}^{l-2r}\dfrac{2^{s}}{s!(l-2r-s)!}\dfrac{I_{s+r+1}(2t)}{t^{r+1}}.
\end{split}\end{equation}

Furthermore, the same procedure and the use of the Hermite polynomials duplication formula \cite{Andrews}

\begin{equation}\label{Hd}
H_{2n}(x,y)=\sum_{r=0}^{n}\binom{n}{r}^{2}r!\;(2y)^{r}\left( H_{n-r}(x,y)\right)^{2}, 
\end{equation}
yields the following identity for Motzkin numbers

\begin{equation}\begin{split}\label{mdupl}
 m_{2n}& =\hat{c} \circ \sum_{r=0}^{n}r!\;\binom{n}{r}^{2}(2\hat{c})^{r} \circ  H_{n-r}(1,\hat{c}) \circ H_{n-r}(1,\hat{c})=\\
& = \sum_{r=0}^{n}\binom{n}{r}^{2}2^{r}r!(n-r)!\sum_{s=0}^{\lfloor\frac{n-r}{2}\rfloor }\dfrac{m_{n-r}^{(r+s+1)}}{(n-r-2s)!s!},
\end{split}\end{equation}
where

\begin{equation}\label{mnmop}
m_{n}^{(q)}=P_{n}^{(q)}(1,1)=\hat{c}^{q} \circ H_{n}(1,\hat{c})
\end{equation}
are associated Motzkin numbers \cite{Blasiak}.\newline
The identification of Motzkin numbers as in eq. \eqref{mnmop}, along with the use of the recurrences of Hermite polynomials, yields, e.g., the identities   

\begin{equation}\begin{split}\label{recmn}
& m_{n+1}^{(q)}=m_{n}^{(q)}+2\;n\;m_{n-1}^{(q+1)},\\
& m_{n+p}=\sum_{s=0}^{\min[n,p]}2^{s}s!\;\binom{p}{s}\;\binom{n}{s}M_{p-s,\;n-s,\;s},\\
& M_{p,\;n,\;t}=p!\sum_{r=0}^{\lfloor\frac{p}{2}\rfloor}\dfrac{m_{n}^{(t+r+1)}}{(p-2r)!r!},
\end{split}\end{equation}
in which, the second identity, has been derived from the Nielsen formula for $H_{n+m}(x,y)$ \cite{Nielsen}.\newline
 
 \section{Final Comments}

In this paper we have shown that a fairly straightforward extension of the formalism put forward in ref. \cite{Blasiak}, allows non trivial progresses in the theory of Motzkin numbers. Further relations can be easily obtained by applying the method we have envisaged as, e.g., 
 
\begin{equation}\label{prmn}
\sum_{s=0}^{n}m_{n-s}\;m_{s}=2\;(n+1)\;m_{n}^{(2)},
\end{equation}
which represents a discrete self-convolution of Motzkin numbers.\newline

 We have also mentioned the existence of the associated Motzkin numbers 

\begin{equation}\label{mnq}
m_{n}^{(q)}=P_{n}^{(q)}(1,1),
\end{equation}
touched on in ref. \cite{Blasiak}. In  the present context they have been introduced on purely algebraic grounds.  Strictly speaking they are not integers and therefore they are not amenable for a combinatorial interpretation however, redefining them as 

\begin{equation}\label{tildem}
\tilde{m}_{n}^{(q)}=\dfrac{(n+q)!}{n!}P_{n}^{(q)}(1,1),
\end{equation}
we obtain for $q=2$ the sequences in OEIS  $(A014531)$, while for $q=3$ the sequences $(A014532)$ and so on.\\

A more appropriate interpretation in combinatorial terms can be obtained by following, e.g., the procedures indicated in ref. \cite{Banderier} and deserves further investigations, out of the scope of the present paper.\newline

We have mentioned in the introduction the  theory of telephone numbers  $T(n)$ \cite{Knuth}, whose importance in chemical Graph theory has been recently emphasized in ref. \cite{Hatz}.  As well known, they can be expressed in terms of ordinary Hermite polynomials, however the use of the two variable extension is more effective. They can indeed be expressed as $T(n)=H_{n}(1,\frac{1}{2})$  .\newline
The use of Hermite polynomials properties, like the index duplication formula, yields

\begin{equation}\label{T2n}
T(2n)=\sum_{r=0}^{n}\;\binom{n}{r}^{2}\;r!\;T(n-r)^{2}.
\end{equation}
 The use of the Hermite numbers $h_{s}$ \cite{Germano} allows the derivation of the following further expression

\begin{equation}\begin{split}\label{hs}
& T(n)=\sum_{s=0}^{n}t_{n,s},\\
& t_{n,s}=\binom{n}{s}\;h_{s}\left( \dfrac{1}{2}\right) ,\\
& h_{s}(y)=y^{\frac{s}{2}}\Gamma\left( \dfrac{s}{2}+2\right)f_{s}=\dfrac{y^{\frac{s}{2}}s!}{\Gamma\left( \dfrac{s}{2}+1\right) }\left| \cos\left( s\;\dfrac{\pi}{2}\right) \right|.  
\end{split}\end{equation}
The coefficients $t_{n,s}$ of the telephone numbers can be arranged in the following triangle, in which, the numbers belonging to the column $s=4$ , $(3, 15, 45, 105, 210, 378,\dots)$,  are identified, in OEIS, with the sequence $A050534$ and the column in $s=6$, $(15, 105, 420,1260, 3150, \dots)$, is just a multiple of $A00910$. \newline

\begin{table}[htbp]
	\centering
	\begin{tabular}{||c | c|c|c|c|c|c|c|c|c||c||} \hline
		\multicolumn{10}{||c||}{\centering $t_{n,s}$ \textbf{coefficients}} & $T{(n)}\; \textbf{telephone numbers}$ \\ \hline
		\toprule
		\multirow{2}{*}{\textbf{Parameter}}  & \multicolumn{9}{c||}{\textbf{s}}&\multirow{2}{*}{$\sum_{s=0}^{n}t_{n,s}$}\\ 
		\cmidrule(lr){3-10}
		&  & \textbf{0} & \textbf{1} & \textbf{2 } & \textbf{3} & \textbf{4} & \textbf{5} & \textbf{6} & \textbf{7}&\\  
		\midrule	
		\multirow{9}{*}{\textbf{n}}  & \textbf{0} & 1 &  &  &  & & & & &\textbf{1}\\ 
		& \textbf{1}  & 1 & 0 &  &  & & &  & & \textbf{1}\\ 
		& \textbf{2} & 1 & 0 & 1 & & & &  & & \textbf{2}\\ 
		& \textbf{3} & 1 & 0 & 3 & 0 &  &  & & & \textbf{4}\\ 
		& \textbf{4} & 1 & 0 & 6 & 0 & 3 &  & & & \textbf{10}\\ 
		& \textbf{5} & 1 & 0 & 10 & 0 & 15 & 0  &  & & \textbf{26}\\ 
		& \textbf{6} & 1 & 0 & 15 & 0 & 45 & 0 & 15 & & \textbf{76}\\ 
		& \textbf{7} & 1 & 0 & 21 & 0 & 105 & 0 & 105 & 0 & \textbf{232}\\ 
		& \dots &  \dots & \dots & \dots    &\dots  &\dots  &\dots    &\dots  &\dots  &\dots\\ \hline
		\bottomrule
	\end{tabular}
	\caption{Telephone Number Coefficients}
	\label{table2}
\end{table}

 The use of the identification with two variable Hermite polynomials opens further perspectives, by exploiting indeed the polynomials (see \cite{Riordan} and references therein)

\begin{equation}\label{Hnm}
H_{n}^{(m)}(x,y)=n!\sum_{r=0}^{\lfloor \frac{n}{m}\rfloor }\dfrac{x^{n-mr}y^{r}}{(n-mr)!r!},
\end{equation}
we can introduce the following generalization of telephone numbers 

\begin{equation}\label{Tnm}
T_{n}^{(m)}=H_{n}^{(m)}\left( 1,\dfrac{1}{m}\right), 
\end{equation} 
with generating function

\begin{equation}\label{gentn}
\sum_{n=0}^{\infty}\dfrac{t^{n}}{n!}T_{n}^{(m)}=e^{t+\frac{1}{m}t^{m}},
\end{equation}
which satisfy the recurrence

\begin{equation}\label{tnmp}
T_{n+1}^{(m)}=T_{n}^{(m)}+\dfrac{n!}{(n-m+1)!}T_{n-m+1}^{(m)}.
\end{equation}
In the case of $m=3$ the numbers $T_{n}^{(3)}= (1, 1, 1, 3, 9, 21, 81, 351, 1233,\dots)$ are identified with OEIS $A001470$, while for $m=4$, the series
$(1, 1, 1, 1, 7, 31, 91, 211, 1681, 12097, \dots)$, corresponds to $A118934$. For $m=5$ the associated series appears to be $A052501$ but should be more appropriately identified with the coefficients of the expansion \eqref{Tnm}, finally the sequence $m=6$ is not reported in OEIS.\newline

 A more accurate analysis of this family of numbers and the relevant interplay with Motzkin will be discussed elsewhere.

\end{document}